\documentclass{amsart} 
\usepackage{latexsym} 
\usepackage{amssymb} 
\usepackage{amsmath} 

\begin{document}


\newtheorem{theorem}{Theorem} 
\newtheorem{problem}{Problem} 
\newtheorem{definition}{Definition} 
\newtheorem{lemma}{Lemma} 
\newtheorem{proposition}{Proposition} 
\newtheorem{corollary}{Corollary} 
\newtheorem{example}{Example} 
\newtheorem{conjecture}{Conjecture} 
\newtheorem{algorithm}{Algorithm} 
\newtheorem{exercise}{Exercise} 
\newtheorem{remarkk}{Remark} 
 
\newcommand{\be}{\begin{equation}} 
\newcommand{\ee}{\end{equation}} 
\newcommand{\bea}{\begin{eqnarray}} 
\newcommand{\eea}{\end{eqnarray}} 
\newcommand{\beq}[1]{\begin{equation}\label{#1}} 
\newcommand{\eeq}{\end{equation}} 
\newcommand{\beqn}[1]{\begin{eqnarray}\label{#1}} 
\newcommand{\eeqn}{\end{eqnarray}} 
\newcommand{\beaa}{\begin{eqnarray*}} 
\newcommand{\eeaa}{\end{eqnarray*}} 
\newcommand{\req}[1]{(\ref{#1})} 
 
\newcommand{\lip}{\langle} 
\newcommand{\rip}{\rangle} 

\newcommand{\uu}{\underline} 
\newcommand{\oo}{\overline} 
\newcommand{\La}{\Lambda} 
\newcommand{\la}{\lambda} 
\newcommand{\eps}{\varepsilon} 
\newcommand{\om}{\omega} 
\newcommand{\Om}{\Omega} 
\newcommand{\ga}{\gamma} 
\newcommand{\rrr}{{\Bigr)}} 
\newcommand{\qqq}{{\Bigl\|}} 
 
\newcommand{\dint}{\displaystyle\int} 
\newcommand{\dsum}{\displaystyle\sum} 
\newcommand{\dfr}{\displaystyle\frac} 
\newcommand{\bige}{\mbox{\Large\it e}} 
\newcommand{\integers}{{\Bbb Z}} 
\newcommand{\rationals}{{\Bbb Q}} 
\newcommand{\reals}{{\rm I\!R}} 
\newcommand{\realsd}{\reals^d} 
\newcommand{\realsn}{\reals^n} 
\newcommand{\NN}{{\rm I\!N}} 
\newcommand{\DD}{{\rm I\!D}} 
\newcommand{\degree}{{\scriptscriptstyle \circ }} 
\newcommand{\dfn}{\stackrel{\triangle}{=}} 
\def\complex{\mathop{\raise .45ex\hbox{${\bf\scriptstyle{|}}$} 
     \kern -0.40em {\rm \textstyle{C}}}\nolimits} 
\def\hilbert{\mathop{\raise .21ex\hbox{$\bigcirc$}}\kern -1.005em {\rm\textstyle{H}}} 
\newcommand{\RAISE}{{\:\raisebox{.6ex}{$\scriptstyle{>}$}\raisebox{-.3ex} 
           {$\scriptstyle{\!\!\!\!\!<}\:$}}} 
 
\newcommand{\hh}{{\:\raisebox{1.8ex}{$\scriptstyle{\degree}$}\raisebox{.0ex} 
           {$\textstyle{\!\!\!\! H}$}}} 

\newcommand{\OO}{\won} 
\newcommand{\calA}{{\mathcal A}} 
\newcommand{\calB}{{\mathcal B}} 
\newcommand{\calC}{{\cal C}} 
\newcommand{\calD}{{\cal D}} 
\newcommand{\calE}{{\mathcal E}} 
\newcommand{\calF}{{\mathcal F}} 
\newcommand{\calG}{{\cal G}} 
\newcommand{\calH}{{\cal H}} 
\newcommand{\calK}{{\cal K}} 
\newcommand{\calL}{{\mathcal L}} 
\newcommand{\calM}{{\cal M}} 
\newcommand{\calO}{{\cal O}} 
\newcommand{\calP}{{\cal P}} 
\newcommand{\calU}{{\mathcal U}} 
\newcommand{\calX}{{\cal X}} 
\newcommand{\calXX}{{\cal X\mbox{\raisebox{.3ex}{$\!\!\!\!\!-$}}}} 
\newcommand{\calXXX}{{\cal X\!\!\!\!\!-}} 
\newcommand{\gi}{{\raisebox{.0ex}{$\scriptscriptstyle{\cal X}$} 
\raisebox{.1ex} {$\scriptstyle{\!\!\!\!-}\:$}}} 
\newcommand{\intsim}{\int_0^1\!\!\!\!\!\!\!\!\!\sim} 
\newcommand{\intsimt}{\int_0^t\!\!\!\!\!\!\!\!\!\sim} 
\newcommand{\pp}{{\partial}} 
\newcommand{\al}{{\alpha}} 
\newcommand{\sB}{{\cal B}} 
\newcommand{\sL}{{\cal L}} 
\newcommand{\sF}{{\cal F}} 
\newcommand{\sE}{{\cal E}} 
\newcommand{\sX}{{\cal X}} 
\newcommand{\R}{{\rm I\!R}} 
\renewcommand{\L}{{\rm I\!L}} 
\newcommand{\vp}{\varphi} 
\newcommand{\N}{{\rm I\!N}} 
\def\ooo{\lip} 
\def\ccc{\rip} 
\newcommand{\ot}{\hat\otimes} 
\newcommand{\rP}{{\Bbb P}} 
\newcommand{\bfcdot}{{\mbox{\boldmath$\cdot$}}} 
 
\renewcommand{\varrho}{{\ell}} 
\newcommand{\dett}{{\textstyle{\det_2}}} 
\newcommand{\sign}{{\mbox{\rm sign}}} 
\newcommand{\TE}{{\rm TE}} 
\newcommand{\TA}{{\rm TA}} 
\newcommand{\E}{{\rm E\,}} 
\newcommand{\won}{{\mbox{\bf 1}}} 
\newcommand{\Lebn}{{\rm Leb}_n} 
\newcommand{\Prob}{{\rm Prob\,}} 
\newcommand{\sinc}{{\rm sinc\,}} 
\newcommand{\ctg}{{\rm ctg\,}} 
\newcommand{\loc}{{\rm loc}} 
\newcommand{\trace}{{\,\,\rm trace\,\,}} 
\newcommand{\Dom}{{\rm Dom}} 
\newcommand{\ifff}{\mbox{\ if and only if\ }} 
\newcommand{\nproof}{\noindent {\bf Proof:\ }} 
\newcommand{\remark}{\noindent {\bf Remark:\ }} 
\newcommand{\remarks}{\noindent {\bf Remarks:\ }} 
\newcommand{\note}{\noindent {\bf Note:\ }}

\newcommand{\boldx}{{\bf x}} 
\newcommand{\boldX}{{\bf X}} 
\newcommand{\boldy}{{\bf y}} 
\newcommand{\boldR}{{\bf R}} 
\newcommand{\uux}{\uu{x}} 
\newcommand{\uuY}{\uu{Y}} 
 
\newcommand{\limn}{\lim_{n \rightarrow \infty}} 
\newcommand{\limN}{\lim_{N \rightarrow \infty}} 
\newcommand{\limr}{\lim_{r \rightarrow \infty}} 
\newcommand{\limd}{\lim_{\delta \rightarrow \infty}} 
\newcommand{\limM}{\lim_{M \rightarrow \infty}} 
\newcommand{\limsupn}{\limsup_{n \rightarrow \infty}} 
 
\newcommand{\ra}{ \rightarrow }

\newcommand{\ARROW}[1] 
  {\begin{array}[t]{c}  \longrightarrow \\[-0.2cm] \textstyle{#1} \end{array} } 
 
\newcommand{\AR} 
 {\begin{array}[t]{c} 
  \longrightarrow \\[-0.3cm] 
  \scriptstyle {n\rightarrow \infty} 
  \end{array}} 
 
\newcommand{\pile}[2] 
  {\left( \begin{array}{c}  {#1}\\[-0.2cm] {#2} \end{array} \right) } 
 
\newcommand{\floor}[1]{\left\lfloor #1 \right\rfloor} 
 
\newcommand{\mmbox}[1]{\mbox{\scriptsize{#1}}} 
 
\newcommand{\ffrac}[2] 
  {\left( \frac{#1}{#2} \right)} 
 
\newcommand{\one}{\frac{1}{n}\:} 
\newcommand{\half}{\frac{1}{2}\:} 
 
\def\le{\leq} 
\def\ge{\geq} 
\def\lt{<} 
\def\gt{>} 
 
\def\squarebox#1{\hbox to #1{\hfill\vbox to #1{\vfill}}} 
\newcommand{\nqed}{\hspace*{\fill} 
           \vbox{\hrule\hbox{\vrule\squarebox{.667em}\vrule}\hrule}\bigskip} 
 
\title{An Anticipating Tangent to a Semimartingale }

\author{ A. S. \"Ust\"unel} 
\maketitle 
\noindent 
{\bf Abstract:}{\small{ We construct a perturbation of identity type
    mapping on an abstract Wiener space where the Cameron-Martin space
    has an orthonormal basis indexed by the jumps a one dimensional
    semimartingale. We then derive a change of variables formula and a
    degree type result for this map.
}}\\ 

\vspace{0.5cm} 

\section{Introduction} 
Denote by $(\Om,\calF,(\calF_t),P)$ a general probability space
satisfying the usual conditions of the theory of stochastic
processes. Let $Z=(Z_t,t\geq 0)$ be a real-valued semimartingale,
denote its Dol\'eans-Dade exponential as $\calE_t(Z)$, i.e., 
$$
\calE_t(Z)=\exp\left(Z_t-\frac{1}{2}\langle Z^c,Z^c\rangle_t\right)\prod_{s\leq
  t}(1+\Delta Z_s) e^{-\Delta Z_s}\,,
$$
where $Z^c$ denotes the continuous local martingale part, $\langle
Z^c,Z^c\rangle_\cdot$ is the associated Doob-Meyer process and $\Delta
Z$ denotes the jumps of $Z$. Since
$P$-almost surely 
\begin{equation}
\label{jump}
\sum_{s\leq t}(\Delta Z_s)^2<\infty
\end{equation}
for any $t\geq 0$, the product term is well-defined. Consequently, for
$P$-almost all $z\in\Om$, the sequence $\{\Delta
Z_s(z):\,s\in[0,t],\,\Delta Z_s(z)\neq 0\}$ can be regarded as the
spectrum of a Hilbert-Schmidt operator and then  the term
$$
D_t=\prod_{s\leq t}(1+\Delta Z_s(z)) e^{-\Delta Z_s(z)}
$$
becomes its modified Carleman-Fredholm determinant. Starting from this
observation, we construct, on a Wiener space an  anticipating perturbation of
identity whose Gaussian Jacobian is the product of $D_t$ with  Wick-Girsanov exponential of this perturbation
of identity. This paper is devoted to studying the properties of this
Jacobian.

\section{Preliminaries}
Let $(W,H,\mu)$ be any Wiener space, denote by $(e_i,i\geq 1)$ a CONB
in the Cameron-Martin space $H$, constructed from the elements of
$W^\star$, i.e., the continuous dual of the Fr\'echet space $W$. Let
$(T_n, n\geq 1)$ be a sequence of stopping times absorbing the jumps
of the semimartingale $Z$ with disjoint graphs, cf. \cite{Del}. We shall construct a
Hilbert-Schmidt operator $\partial u_t$, whose eigenvalues will consist of the jumps
of $Z$: define
$$
\partial u_t=\sum_{T_n\leq t}\Delta Z_{T_n}\,e_n\otimes e_n\,.
$$
Because of the property (\ref{jump}), $\partial u_t$ is a measurable
map defined on $\R_+\times\Om$ with values in $H\tilde{\otimes}_2H$,
where the latter denotes the completed Hilbert-Schmidt tensor
product. In fact
$$
\| \partial u_t(z)\|_2^2=\sum_{T_n(z)\leq t}(\Delta Z_{T_n}(z))^2 \leq \sum_{s\leq t}(\Delta Z_s(z))^2<\infty
$$
$P$-almost surely. Since this operator is also symmetric, for fixed
$z$, we can construct an $H$-valued map, whose Sobolev derivative on
the Wiener space will be equal to $\partial u_t$ explicitly. In fact,
it suffices to put
$$
u_t(w,z)=\sum_{T_n(z)\leq t}\Delta Z_{T_n}(z)\delta e_n(w)\,e_n\,,
$$
where $\delta e_n$ denotes the Gaussian divergence of $e_n\in H$. Note
that, defining 
$$
u^\eps_t(w,z)\sum_{T_n\leq t,|\Delta Z_{T_n}|>\eps}\Delta
Z_{T_n}(z)\delta e_n(w)\,e_n
$$
we get well-defined,  $H$-valued, measurable maps. For $\eps,\eta>0$,
we have 
$$
|u^\eps_t-u^\eta_t|_H^2=|u^\eps_t|_H^2+|u^\eta_t|_H^2-2\sum_{T_n(z)\leq
  t,\,\eps\wedge \eta<|\Delta Z_{T_n}|\leq\eps\vee \eta}(\Delta Z_{T_n}(z))^2(\delta
e_n(w))^2
$$
since the graphs of $(T_n)$ are disjoint. Consequently, for $P$-almost
all $z$, $\lim_{\eps\to 0} u^\eps_t(\cdot,z)=u_t(\cdot,z)$ in
$L^2(\mu,H)$, since $\lim\nabla u^\eps_t(\cdot,z)$ exists also,
$u_t(\cdot,z)\in \DD_{p,1}(H)$ for any $p>1$. We have also a more
delicate result:
\begin{theorem}
\label{main-thm}
The mapping $U_t(w,z)=w+u_t(w,z)$ defined on $W\times \Om$ with values
in $W$, where $u_t$ is given as
$$
u_t(w,z)=\sum_{T_n(z)\leq t}\Delta Z_{T_n}\delta e_n(w)e_n\,,
$$
is $\mu\times P$-almost surely well-defined perturbation of identity
on the Wiener space $W$, besides, for $P$-almost all $z$, the partial
map 
$$
w\to u_t(w,z)
$$
is an $H-C^1$-map on $W$, i.e., the map $h\to u_t(w+h,z)$ is a
$C^1$-map on $H$ outside a set of capacity zero (even an $H$-invariant
set).
\end{theorem}
\nproof
For any $h\in H$, we have 
$$
u^\eps_t(w+h,z)=u^\eps_t(w+,z)+\nabla u^\eps_t(z)h
$$
$P$-a.s. for any $w\in W$ (recall that $e_n\in W^\star$). Hence
\beaa
|u^\eps_t(w+h,z)-u^\eta_t(w+h,z)|_H^2&\leq&
2|u^\eps_t(w,z)-u^\eta_t(w,z)|_H^2\\
&&+2|h|_H^2\|\nabla u^\eps_t(z)-\nabla u^\eta_t(z)\|_2^2\,.
\eeaa
Hence $(h\to u^\eps_t(w+h,z),\eps>0)$ converges to $h\to
u_t(w+h,z)$ uniformly w.r.to $h\in H$ and this implies that $h\to
u_t(w+h,z)$ is an $H-C$-map for almost all $z\in\Om$. Moreover,
$\nabla u_t(w,z)$ does not depend on $w$, hence $h\to u_t(w,z)$ is
even $H-C^\infty$.
\nqed

\noindent
In the sequel $\La_t=\La_t(w,z)$ will represent the Gaussian Jacobian
associated to the perturbation of identity $U_t(w,z)$, namely
\beaa
\La_t(w,z)&=&\dett(I_H+\nabla u_t(w,z))\exp(-\delta
u_t(w,z)-\half|u_t(w,z)|_H^2)\\
&=&\prod_{T_n(z)\leq t}(1+\Delta Z_{T_n})e^{-\Delta
  Z_{T_n}}\,\exp\left(-\sum_{T_n\leq t}\Delta Z_{T_n}((\delta
  e_n)^2-1)-\half\sum_{T_n\leq t}(\Delta Z_{T_n})^2(\delta e_n)^2\right)
\eeaa

\noindent
We have the following degree-type result (cf. \cite{A-E} or \cite{BOOK}):
\begin{theorem}
\label{Sard-thm}
Let us denote by $M_t$ the set of nondegeneracy of $\La_t$, namely 
$$
M_t=\bigcap_{T_n\leq t}\{z\in \Om: \,\Delta Z_{T_n}\neq -1\}\,.
$$
Then we have 
$$
\mu(U_t(M_t^c))=0
$$
and
\begin{equation}
\label{chg-of-v}
\int f(z,U_t(w,z))|\La_t(w,z)|g(w,z)d\mu dP=\int f(w,z)\left(\sum_{y\in
    U_t^{-1}(\cdot,z)\{w\}}g(y,z)\right)d\mu dP
\end{equation}
for any positive, measurable functions on $W\times\Om$. Consequently,
we have 
\begin{equation}
\label{deg-eqn}
\int_W f\circ U_t(w,z)\,\La_t(w,z)d\mu=\int_W fd\mu\int_W \La_td\mu
\end{equation}
$P$-almost surely and hence 
\begin{equation}
\label{deg-eqn1}
\int_W \La_td\mu=\sum_{y\in  U_t^{-1}(\cdot,z)\{w\}}\sign\La_t(y,z)=|\{ U_t^{-1}(\cdot,z)\{w\}\}|\sign(D_t(z))
\end{equation}
$P\times\mu$-almost surely, where $|\{\ldots\}|$ denotes the
cardinality of the set inside the brackets.
Moreover, in the case of classical Wiener space,  the mapping
$w\to U_t(w,z)$  is $\mu$-almost surely invertible on the set $\{z\in\Om:\,E_\mu[|\La_t|]=1\}$, and its
(measurable) inverse, denoted as $(\tau,w)\to V_t(w,z)(\tau)$, is the
strong solution of the following functional stochastic differential equation
$$
V_t(w,z)(\tau)=W_\tau-\int_0^\tau\dot{u}_t(V_t(w,z),z)(s)ds\,.
$$

\end{theorem}
\nproof
Since $U_t$ is an $H-C^1$-map $P$-almost surely, the Sard inequality
implies that  the image  under
itself of  its set of degeneracy has zero $\mu$-measure  and the
change of variables formula (\ref{chg-of-v}) holds
cf. \cite{BOOK}, Chapter 4, Section 4.4. The relation (\ref{deg-eqn})
follows from Theorem 9.2.3 of \cite{BOOK} and finally the equality
(\ref{deg-eqn1}) is a sequence of the relations (\ref{deg-eqn}) and
(\ref{deg-eqn1}) . As explained in \cite{INV}, the equivalence of the
measures  $U_t\mu$ and $\mu$, combined with the fact that the
cardinality of $U_t^{-1}(\cdot,z)\{w\}$ being equal to one
$\mu$-almost everywhere  implies the
almost sure invertibility the map $w\to U_t(w,z)$ $P$-almost
everywhere and the functional stochastic differential equation is a
direct consequence of this fact.
\nqed\newline
\remark Note that the delicate point in the above result comes from
the lack of information about the instants of the  distorted noise
posterior to $t\in \R_+$.

\noindent
We can view $(\La_t\in[0,1])$ as a semimartingale with respect to the
filtration $(\calF_t)$ . First let us define the process
$(e(t),t\in[0,1])$ as
$$
e(t)=\sum_n e_n\,1_{[T_n,T_{n+1}[}(t)\,.
$$
We have
\begin{theorem}
\label{sde-thm}
Assume that $Z$ is a pure jump semimartingale, then the process
$(t,z)\to \La_t(w,z)$ satisfies the following stochastic differential
equation:
\beaa
\La_t&=&\La_0-\int_0^t \La_{s-}\left[(\delta e_s)^2(1+\half\Delta
  Z_s)-1\right]dZ_s\\
&&+\sum_{s\leq t}\La_{s-}e^{\Delta Z_s[(\delta e_s)^2(1-\half\Delta
  Z_s)-1]}(e^{-\Delta Z_s} (1+\Delta Z_s)-1)
\eeaa
\end{theorem}

\vspace{2cm}
\footnotesize{
\noindent
A. S. \"Ust\"unel,\\
Mathematics Departement, Bilkent University, 
Ankara, Turkey,\\
ustunel@fen.bilkent.edu.tr}

\end{document}